\documentclass[11pt, a4paper]{article}
\usepackage[cp1251]{inputenc}
\usepackage{amsmath} \usepackage{euscript} \usepackage{marvosym}
\usepackage[dvipsnames]{xcolor}
\usepackage{mymatrx6}%\let\mymatrixII\mymatrix
\usepackage{comment}
\usepackage{graphicx}
\usepackage{longtable}
\usepackage{mathrsfs}
\usepackage{amssymb, amscd}
\usepackage{adigraph}
\usepackage{nicefrac}
\def\mymatrix{\MyMatrixwithdelims..}
\Linestrue\Autonumtrue
\usepackage{mathrsfs}
\def\bluecommand{\textcolor[cmyk]{1,0,0,0}}
\def\redcommand{\textcolor[cmyk]{0,1,1,0}}
\def\greencommand{\textcolor[cmyk]{1,0,1,0}}
\def\yellowcommand{\textcolor[cmyk]{0,0.25,1,0}}
\def\graycommand{\textcolor[cmyk]{0,0,0,0.2}}

\newtheorem{theorem}{Theorem.}
\newtheorem{lemma}[theorem]{Lemma.}

\begin{document}

\newcounter{bnomer} \newcounter{snomer}
\newcounter{bsnomer}
\setcounter{bnomer}{0}
\renewcommand{\thesnomer}{\thebnomer.\arabic{snomer}}
\renewcommand{\thebsnomer}{\thebnomer.\arabic{bsnomer}}
\renewcommand{\refname}{\begin{center}\large{\textbf{References}}\end{center}}

\setcounter{MaxMatrixCols}{14}

\newcommand\restr[2]{{% we make the whole thing an ordinary symbol
  \left.\kern-\nulldelimiterspace % automatically resize the bar with \right
  #1 % the function
  %\vphantom{\big|} % pretend it's a little taller at normal size
  \right|_{#2} % this is the delimiter
}}

\newcommand{\sect}[1]{%
\setcounter{snomer}{0}\setcounter{bsnomer}{0}
\refstepcounter{bnomer}
\par\bigskip\begin{center}\large{\textbf{\arabic{bnomer}. {#1}}}\end{center}}
\newcommand{\sst}[1]{%
\refstepcounter{bsnomer}
\par\bigskip\textbf{\arabic{bnomer}.\arabic{bsnomer}. {#1}}\par}
\newcommand{\defi}[1]{%
\refstepcounter{snomer}
\par\medskip\textbf{Definition \arabic{bnomer}.\arabic{snomer}. }{#1}\par\medskip}
\newcommand{\theo}[2]{%
\refstepcounter{snomer}
\par\textbf{Theorem \arabic{bnomer}.\arabic{snomer}. }{#2} {\emph{#1}}\hspace{\fill}$\square$\par}
\newcommand{\theobbp}[1]{%
\refstepcounter{snomer}
\par\textbf{Theorem \arabic{bnomer}.\arabic{snomer}.} {\emph{#1}}\hspace{\fill}$\square$\par}
\newcommand{\mtheop}[2]{%
\refstepcounter{snomer}
\par\textbf{Theorem \arabic{bnomer}.\arabic{snomer}. }{\emph{#1}}
\par\textsc{Proof}. {#2}\hspace{\fill}$\square$\par}
\newcommand{\mcorop}[2]{%
\refstepcounter{snomer}
\par\textbf{Corollary \arabic{bnomer}.\arabic{snomer}. }{\emph{#1}}
\par\textsc{Proof}. {#2}\hspace{\fill}$\square$\par}
\newcommand{\mtheo}[1]{%
\refstepcounter{snomer}
\par\medskip\textbf{Theorem \arabic{bnomer}.\arabic{snomer}. }{\emph{#1}}\par\medskip}
\newcommand{\theobn}[1]{%
\par\medskip\textbf{Theorem. }{\emph{#1}}\par\medskip}
\newcommand{\theoc}[2]{%
\refstepcounter{snomer}
\par\medskip\textbf{Theorem \arabic{bnomer}.\arabic{snomer}. }{#1} {\emph{#2}}\par\medskip}
\newcommand{\mlemm}[1]{%
\refstepcounter{snomer}
\par\medskip\textbf{Lemma \arabic{bnomer}.\arabic{snomer}. }{\emph{#1}}\par\medskip}
\newcommand{\mprop}[1]{%
\refstepcounter{snomer}
\par\medskip\textbf{Proposition \arabic{bnomer}.\arabic{snomer}. }{\emph{#1}}\par\medskip}
\newcommand{\theobp}[2]{%
\refstepcounter{snomer}
\par\textbf{Theorem \arabic{bnomer}.\arabic{snomer}. }{#2} {\emph{#1}}\par}
\newcommand{\theop}[2]{%
\refstepcounter{snomer}
\par\textbf{Theorem \arabic{bnomer}.\arabic{snomer}. }{\emph{#1}}
\par\textsc{Proof}. {#2}\hspace{\fill}$\square$\par}
\newcommand{\theosp}[2]{%
\refstepcounter{snomer}
\par\textbf{Theorem \arabic{bnomer}.\arabic{snomer}. }{\emph{#1}}
\par\textsc{Sketch of the proof}. {#2}\hspace{\fill}$\square$\par}
\newcommand{\exam}[1]{%
\refstepcounter{snomer}
\par\medskip\textbf{Example \arabic{bnomer}.\arabic{snomer}. }{#1}\par\medskip}
\newcommand{\deno}[1]{%
\refstepcounter{snomer}
\par\textbf{Notation \arabic{bnomer}.\arabic{snomer}. }{#1}\par}
\newcommand{\lemm}[1]{%
\refstepcounter{snomer}
\par\textbf{Lemma \arabic{bnomer}.\arabic{snomer}. }{\emph{#1}}\hspace{\fill}$\square$\par}
\newcommand{\lemmp}[2]{%
\refstepcounter{snomer}
\par\medskip\textbf{Lemma \arabic{bnomer}.\arabic{snomer}. }{\emph{#1}}
\par\textsc{Proof}. {#2}\hspace{\fill}$\square$\par\medskip}
\newcommand{\coro}[1]{%
\refstepcounter{snomer}
\par\textbf{Corollary \arabic{bnomer}.\arabic{snomer}. }{\emph{#1}}\hspace{\fill}$\square$\par}
\newcommand{\mcoro}[1]{%
\refstepcounter{snomer}
\par\textbf{Corollary \arabic{bnomer}.\arabic{snomer}. }{\emph{#1}}\par\medskip}
\newcommand{\corop}[2]{%
\refstepcounter{snomer}
\par\textbf{Corollary \arabic{bnomer}.\arabic{snomer}. }{\emph{#1}}
\par\textsc{Proof}. {#2}\hspace{\fill}$\square$\par}
\newcommand{\nota}[1]{%
\refstepcounter{snomer}
\par\medskip\textbf{Remark \arabic{bnomer}.\arabic{snomer}. }{#1}\par\medskip}
\newcommand{\propp}[2]{%
\refstepcounter{snomer}
\par\medskip\textbf{Proposition \arabic{bnomer}.\arabic{snomer}. }{\emph{#1}}
\par\textsc{Proof}. {#2}\hspace{\fill}$\square$\par\medskip}
\newcommand{\hypo}[1]{%
\refstepcounter{snomer}
\par\medskip\textbf{Conjecture \arabic{bnomer}.\arabic{snomer}. }{\emph{#1}}\par\medskip}
\newcommand{\prop}[1]{%
\refstepcounter{snomer}
\par\textbf{Proposition \arabic{bnomer}.\arabic{snomer}. }{\emph{#1}}\hspace{\fill}$\square$\par}

\newcommand{\proof}[2]{%
\par\medskip\textsc{Proof{#1}}. \hspace{-0.2cm}{#2}\hspace{\fill}$\square$\par\medskip}

\makeatletter
\def\iddots{\mathinner{\mkern1mu\raise\p@
\vbox{\kern7\p@\hbox{.}}\mkern2mu
\raise4\p@\hbox{.}\mkern2mu\raise7\p@\hbox{.}\mkern1mu}}
\makeatother

\newcommand{\okr}[2]{%
\refstepcounter{snomer}
\par\medskip\textbf{{#1} \arabic{bnomer}.\arabic{snomer}. }{\emph{#2}}\par\medskip}

\newcommand{\Ind}[3]{%
\mathrm{Ind}_{#1}^{#2}{#3}}
\newcommand{\Res}[3]{%
\mathrm{Res}_{#1}^{#2}{#3}}
\newcommand{\epsi}{\varepsilon}
\newcommand{\tri}{\triangleleft}
\newcommand{\Supp}[1]{%
\mathrm{Supp}(#1)}
\newcommand{\NSupp}[1]{%
\mathrm{NSupp}(#1)}
\newcommand{\SSu}[1]{%
\mathrm{SingSupp}(#1)}

\newcommand{\lee}{\leqslant}
\newcommand{\gee}{\geqslant}
\newcommand{\reg}{\mathrm{reg}}
\newcommand{\Dyn}{\mathrm{Dyn}}
\newcommand{\Ann}{\mathrm{Ann}\,}
\newcommand{\Cent}[1]{\mathbin\mathrm{Cent}({#1})}
\newcommand{\PCent}[1]{\mathbin\mathrm{PCent}({#1})}
\newcommand{\Irr}[1]{\mathbin\mathrm{Irr}({#1})}
\newcommand{\Exp}[1]{\mathbin\mathrm{Exp}({#1})}
\newcommand{\empr}[2]{[-{#1},{#1}]\times[-{#2},{#2}]}
\newcommand{\sreg}{\mathrm{sreg}}
\newcommand{\ilm}{\varinjlim}
\newcommand{\wdth}{\mathrm{wd}}
\newcommand{\plm}{\varprojlim}
\newcommand{\codim}{\mathrm{codim}\,}
\newcommand{\GKdim}{\mathrm{GKdim}\,}
\newcommand{\chara}{\mathrm{char}\,}
\newcommand{\rk}{\mathrm{rk}\,}
\newcommand{\chr}{\mathrm{ch}\,}
\newcommand{\Ker}{\mathrm{Ker}\,}
\newcommand{\id}{\mathrm{id}}
\newcommand{\Ad}{\mathrm{Ad}}
\newcommand{\Gh}{\mathrm{Gh}}
\newcommand{\col}{\mathrm{col}}
\newcommand{\row}{\mathrm{row}}
\newcommand{\high}{\mathrm{high}}
\newcommand{\low}{\mathrm{low}}
\newcommand{\pho}{\hphantom{\quad}\vphantom{\mid}}
\newcommand{\fho}[1]{\vphantom{\mid}\setbox0\hbox{00}\hbox to \wd0{\hss\ensuremath{#1}\hss}}
\newcommand{\wt}{\widetilde}
\newcommand{\wh}{\widehat}
\newcommand{\ad}[1]{\mathrm{ad}_{#1}}
\newcommand{\tr}{\mathrm{tr}\,}
\newcommand{\GL}{\mathrm{GL}}
\newcommand{\SL}{\mathrm{SL}}
\newcommand{\SO}{\mathrm{SO}}
\newcommand{\Or}{\mathrm{O}}
\newcommand{\Sp}{\mathrm{Sp}}
\newcommand{\SuppD}{\mathbb{S}\mathrm{upp}}
\newcommand{\Sa}{\mathrm{S}}
\newcommand{\Sing}{\mathrm{Sing}}
\newcommand{\Ua}{\mathrm{U}}
\newcommand{\Andre}{\mathrm{Andre}}
\newcommand{\Aord}{\mathrm{Aord}}
\newcommand{\Mat}{\mathrm{Mat}}
\newcommand{\Stab}{\mathrm{Stab}}
\newcommand{\htt}{\mathfrak{h}}
\newcommand{\Hei}{\mathrm{Hei}}
\newcommand{\spt}{\mathfrak{sp}}
\newcommand{\slt}{\mathfrak{sl}}
\newcommand{\sot}{\mathfrak{so}}

\newcommand{\vfi}{\varphi}
\newcommand{\aad}{\mathrm{ad}}
\newcommand{\vpi}{\varpi}
\newcommand{\teta}{\vartheta}
\newcommand{\Bfi}{\Phi}
\newcommand{\Fp}{\mathbb{F}}
\newcommand{\Rp}{\mathbb{R}}
\newcommand{\Zp}{\mathbb{Z}}
\newcommand{\Cp}{\mathbb{C}}
\newcommand{\Ap}{\mathbb{A}}
\newcommand{\Pp}{\mathbb{P}}
\newcommand{\Kp}{\mathbb{K}}
\newcommand{\Np}{\mathbb{N}}
\newcommand{\ut}{\mathfrak{u}}
\newcommand{\at}{\mathfrak{a}}
\newcommand{\glt}{\mathfrak{gl}}
\newcommand{\hei}{\mathfrak{hei}}
\newcommand{\nt}{\mathfrak{n}}
\newcommand{\kt}{\mathfrak{k}}
\newcommand{\mt}{\mathfrak{m}}
\newcommand{\rt}{\mathfrak{r}}
\newcommand{\rad}{\mathfrak{rad}}
\newcommand{\bt}{\mathfrak{b}}
\newcommand{\unt}{\underline{\mathfrak{n}}}
\newcommand{\gt}{\mathfrak{g}}
\newcommand{\vt}{\mathfrak{v}}
\newcommand{\pt}{\mathfrak{p}}
\newcommand{\Xt}{\mathfrak{X}}
\newcommand{\Po}{\mathcal{P}}
\newcommand{\PV}{\mathcal{PV}}
\newcommand{\Uo}{\EuScript{U}}
\newcommand{\Fo}{\EuScript{F}}
\newcommand{\Do}{\EuScript{D}}
\newcommand{\Eo}{\EuScript{E}}
\newcommand{\Jo}{\EuScript{J}}
\newcommand{\Iu}{\mathcal{I}}
\newcommand{\Mo}{\mathcal{M}}
\newcommand{\Nu}{\mathcal{N}}
\newcommand{\Ro}{\mathcal{R}}
\newcommand{\Co}{\mathcal{C}}
\newcommand{\Ko}{\mathcal{K}}
\newcommand{\So}{\mathcal{S}}
\newcommand{\Lo}{\mathcal{L}}
\newcommand{\Ou}{\mathcal{O}}
\newcommand{\Uu}{\mathcal{U}}
\newcommand{\Tu}{\mathcal{T}}
\newcommand{\Au}{\mathcal{A}}
\newcommand{\Vu}{\mathcal{V}}
\newcommand{\Du}{\mathcal{D}}
\newcommand{\Bu}{\mathcal{B}}
\newcommand{\Sy}{\mathcal{Z}}
\newcommand{\Sb}{\mathcal{F}}
\newcommand{\Gr}{\mathcal{G}}
\newcommand{\Xu}{\mathcal{X}}
\newcommand{\Op}{\mathbb{O}}
\newcommand{\chv}{\mathrm{chv}}
\newcommand{\rtc}[1]{C_{#1}^{\mathrm{red}}}

\newcommand{\JSpec}[1]{\mathrm{JSpec}\,{#1}}
\newcommand{\MSpec}[1]{\mathrm{MSpec}\,{#1}}
\newcommand{\PSpec}[1]{\mathrm{PSpec}\,{#1}}
\newcommand{\APbr}[1]{\mathrm{span}\{#1\}}
\newcommand{\APbre}[1]{\langle #1\rangle}
\newcommand{\APro}[1]{\setcounter{AP}{#1}\Roman{AP}}\newcommand{\apro}[1]{{\rm\setcounter{AP}{#1}\roman{AP}}}
\newcommand{\ot}{\xleftarrow[]{}}
\newcounter{AP}

%\contentsname{Table of contents}

\author{Mikhail Ignatev\and Leonid Titov}
\date{}
\title{On the number of irreducible representations\\ for finite unipotent Heisenberg-type groups}\maketitle
\begin{abstract} Let $U$ be an algebraic subgroup of the group of $n\times n$ upper-triangular matrices with units on the diagonal over a finite field of large enough characteristic, and $\nt$ be the Lie algebra of~$U$. The main tool in representation theory of $U$ is the orbit method, which classifies irreducible representations of the group $U$ in terms of coadjoint orbits on the dual space $\nt^*$. We consider two types of generalizations of the Heisenberg group, namely, generalized Heisenberg groups defined with an arbitrary bilinear form, and certain subgroups in maximal unipotent subgroups of classical orthogonal algebraic groups. We provide a way to calculate the number of irreducible representations of such groups. It turned out that this number is a polynomial in $q-1$ with nonnegative integer coefficients, which agrees with Isaacs' conjecture.

\end{abstract}
%\tableofcontents\newpage
%\maketitle

%\tableofcontents\newpage

\vspace{0.2cm}

\sect{Introduction}\label{sect:intro}\addcontentsline{toc}{subsection}{\ref{sect:intro}. Introduction}

Throughout the paper, $U$ denotes an algebraic subgroup of the group $U_n(q)$ of upper-triangular $n\times n$ matrices over a finite field $\Fp_q$ of sufficiently large characteristic $p$. The group $U$ is a unipotent affine group and a $q$-power degree group, which means that every irreducible complex finite-dimensional representation of $U$ is of dimension a power of $q$ (in fact, the dimension of such a representation is $q^{2e}$ for some $e\geq0$). We will denote the set of irreducible characters of $U$ by $\Irr{U}$. A long-standing Higman's conjecture claims that $|\Irr{U}|$ is a polynomial in $q$ with integer coefficients. Later, Lehrer made a stronger conjecture, which states that $|\Irr{U}_e|$ is a polynomial in $q$ with integer coefficients, where $\Irr{U}_e$ is the set of irreducible characters of degree~$q^{2e}$. Finally, Isaacs formulated an even stronger conjecture: $|\Irr{U}_q|$ is a polynomial in $v=q-1$ with non-negative integer coefficients. (Actually, Higman, Lehrer and Isaacs formulated their conjectures only for $U=U_n(q)$, but, of course, they are very interesting for all finite unipotent affine groups $U$.)

In past years, a significant progress in investigation of these conjectures was been made, see the next section for the detail. In this paper, we prove that Isaacs' conjecture holds for generalized Heisenberg groups (see Section~\ref{sect:gener_Hei} for precise definitions), as well as for so-called classical orthogonal 2-layered Heisenberg groups (see Section~\ref{Sfdcn} for the detailed description). Both of these types of groups generalize the classical Heisenberg group, which plays a central role in various branches of mathematics, such as representation theory, algebraic geometry, quantum mechanics, etc. The main technical tool we use for generalized Heisenberg group is the Kirillov's orbit method which we briefly recall in the next section, while for orthogonal 2-layered Heisenberg groups we use the Mackey little group method for semi-direct products, see Section~\ref{sect:Mackey}. The main results of the paper are formulated in Proposition~\ref{prop:gen_Hei_orbits}, Corollary~\ref{coro:gen_Hei_Isaacs} and Theorem~\ref{theo:2_Hei_Isaacs}.We thank Mikhail Venchakov for useful discussions which helped us to improve the text of the paper.

\let\thefootnote\relax\footnote{The work of the first author is an output of a research project implemented as part of the Basic Research Program at HSE University.}

\sect{Isaacs' conjecture\label{sect:number_chars}}\addcontentsline{toc}{subsection}{\ref{sect:number_chars}. Isaacs' conjecture}

As above, let $U$ be an algebraic subgroup of the group $U_n(q)$ (we will assume that the characteristic $p$ of the ground field $\Fp_q$ satisfies the condition $p\geq n$).  For example, one can put $U$ to be a maximal unipotent subgroup in a classical simple group $G$ over $\Fp_q$ (or, equivalently, a Sylow $p$-subgroup of $G$). Let $\nt$ be the Lie algebra of $U$, and $\nt^*$ be the dual space to $\nt$. Clearly, $U$ acts on $\nt$ via the adjoint action:
$$g\in U,~x\in\nt\mapsto gxg^{-1}.$$
The dual action of $N$ on $\nt^*$ is called coadjoint:
$$(g.\lambda)(x)=\lambda(g^{-1}xg),~g\in U,~x\in\nt,~\lambda\in\nt^*.$$
It is well known that a coadjoint orbit $\Omega\subseteq\nt^*$ is an affine variety defined over $\Fp_1$, and its dimension $\dim\Omega$ equals $2e$ for some $e\geq0$ (and $|\Omega|=q^{2e}$, because the orbit $\Omega$ is isomorphic to the affine space of dimension $2e$). 

Fix a nontrivial group homomorphism $\theta\colon\Fp_q\to\Cp^{\times}$, where $\Cp^{\times}$ is the multiplicative group of $\Cp$. Let $\ln\colon U\to\nt$ be the inverse map to the exponential map. Given a coadjoint $U$-orbit $\Omega\subset\nt^*$, put
$$\chi_{\Omega}(g)=\dfrac{1}{\sqrt{|\Omega|}}\sum_{\lambda\in\Omega}\theta(\lambda(\ln g)),~g\in U.$$
The orbit method created by A.A. Kirillov in 1962 \cite{Kirillov62} claims that
the function $\chi_{\Omega}\colon U\to\Cp$ is an irreducible complex character of the group $U$, each irreducible character of $U$ has such a form, and $\chi_{\Omega_1}=\chi_{\Omega_2}$ if and only if $\Omega_1=\Omega_2$. In other words,
the map $\Omega\mapsto\chi_{\Omega}$ establishes a bijection between the set of coadjoint $U$-orbits on $\nt^*$ and the set $\Irr{U}$ of irreducible complex characters of the group $U$. It is clear that if $\dim\Omega=2e$ then the \emph{degree} $\deg\chi_{\Omega}$ of $\chi_{\Omega}$ (i.e., the complex dimension of the corresponding irreducible representation) equals $q^e$.

A longstanding G. Higman's conjecture \cite{Higman60} states that the number of conjugacy classes for~$U_n(q)$ is a polynomial in~$q$ (note that $U_n(q)$ is a maximal unipotent subgroup in the simple group $G=\SL_n(\Fp_q)$ with the root system $\Phi=A_{n-1}$). It is easy to check that the number of conjugacy classes for $U$ coincides with the number of coadjoint $U$-orbits, or, equivalently, with the number of irreducible characters of $U$. Of course, it is interesting to consider an analogue of Higman's conjecture for an arbitrary $U$, not only for $U_n(q)$; we will denote the number of irreducible characters of the group $U$ by $O(q)$. Forty years after Higman, G. Lehrer conjectured \cite{Lehrer74} that, for $U_n(q)$, even the number $O_e(q)$ of characters of degree $q^e$ is polynomial in $q$. In 2007, I.M. Isaacs made a stronger conjecture \cite{Isaacs07} that, for $U_n(q)$, $O_e(q)$ is in fact polynomial in $q-1$ with nonnegative integer coefficients.

In past twenty five years, a significant progress in studying these conjectures has been made for $U=U_n(q)$ (we will write $O_{n,e}(q)$ instead of $O_e(q)$ in this case). Original Higman's conjecture was checked for $n\leq13$ in 2003 by A. Vera-Lopez and and J.M. Arregi in \cite{VeraLopezArregi03} by computer computations. In 2007, Isaacs proposed conjectural polynomials for $O_{n,e}(q)$ with $n\leq9$; he also computed explicitly $O_{n,\mu(n)}(q)$ and $O_{n,\mu(n)-1}(q)$, where $q^{\mu(n)}$ is the maximal possible degree of an irreducible character for $U_n(q)$. Note that these formulas were proved in 1997 by M. Marjoram in his unpublished thesis \cite{Marjoram97'},~\cite{Marjoram97}; the formula for $O_{n,\mu(n)-1}(q)$ also follows from A. Panov's classification of $2(\mu(n)-1)$-dimensional orbits presented in \cite{IgnatevPanov09}. In 2010, A. Evseev \cite{Evseev10} computed $O_{n,e}(q)$ for $n\leq13$ confirming Isaacs' formulas (and so Isaacs' conjecture). In 1999, Marjoram calculated $O_{n,1}(q)$ and $O_{n,2}(q)$ \cite{Marjoram99} (clearly, $O_{n,0}(q)=q^{n-1}$), while $O_{n,3}(q)$ was calculated by M. Loukaki in 2011 \cite{Loukaki11}. The results of Marjoram were independently rederived by T. Le in 2010 \cite{Le10}. E. Marberg verified Isaacs' conjecture for $e\leq8$ for arbitrary $n$ in 2011 \cite{Marberg11}.

For Sylow subgroups $U$ in other simple groups $G$ over $\Fp_q$, the situation is as follows. In 1999, Marjoram computed $O_{\mu(\Phi)}(q)$ for orthogonal group (i.e., for $\Phi=B_n$ or $D_n$), where $q^{\mu(\Phi)}$ is the maximal possible degree of an irreducible character of $U$, see \cite{Marjoram99}. For the symplectic case (i.e., for $\Phi=C_n$), the formula for $O_{\mu(\Phi)}(q)$ follows from the results of C.A.M. Andr\`e and A.-M. Neto on so-called supercharacters published in 2006 \cite{AndreNeto06} (see also \cite{Venchakov26}). In 2016, S.M. Goodwin, P. Mosch and G. R\"ohrle \cite{GoodwinMoschRoehrle16} calculated $O_e(q)$ and proved Isaacs' conjecture for all possible $e$ and all finite simple groups $G$ of rank $\leq 8$, except $E_8$.
\exam{Their calculations show, e.g., that, for $\Phi=F_4$, one has $O_4(q)=v^8 + 8v^7 + 28v^6 + 58v^5 + 79v^4 + 66v^3 + 24v^2 + 2v$, where $v=q-1$, while for $\Phi=E_7$, $O_{13}(q)=3v^7 + 24v^6 + 63v^5 + 68v^4 + 28v^3 + 3v^2$.}

We also would like to mention the papers \cite{GoodwinLeMagaardPaolini16} and \cite{LeMagaardPaolini20} of S.M. Goodwin, T. Le, K. Magaard and A. Paolini, where the case $F_4$ is considered without any restrictions on the characteristic (see also \cite{Sur26}), as well as the paper \cite{Szechtman06}, where F. Szechtman constructed certain explicit families of characters for symplectic and orthogonal cases. The case $D_4$ was considered in detail in \cite{GoodwinLeMagaard}.

In the paper \cite{IgnatevPetukhov25}, A. Petukhov together with the first author classified coadjoint orbits of dimension less or equal to 6 for maximal unipotent subgroups in simple groups of types $B_n$, $C_n$ and $D_n$. It follows from this classification that Isaacs' conjecture holds in this case.

\exam{It follows from this classification that, for $\Phi=A_{n-1}$, $n\ge7$, one has
\begin{equation*}
\label{E:exfq}
\begin{split}
%O_1(q)&=(n-2)(v+1)^{n-3}v+(n-3)(v+1)^{n-4}(v^3+v^2)\\
%&=(v+1)^{n-4}(n(v^3+2v^2+v)-(3v^3+5v^2+2v)),\\
O_3(q)&=(n-7)q^{n+1}+(n^2-14n+52)q^n+\dfrac{n^3-30n^2+293n-906}{6}q^{n-1}\\
&-\dfrac{n^3-21n^2+164n-446}{2}q^{n-2}
+\dfrac{n^3-21n^2+150n-370}{2}q^{n-3}\\
&-\dfrac{n^3-27n^2+206n-498}{6}q^{n-4}
-\dfrac{n^2-11n+30}{2}q^{n-5}\\
&=(v+1)^{n-5}\left((n-7)v^6+(n^2-8n+10)v^5+\dfrac{n^3-37n+24}{6}v^4+\dfrac{n^3+3n^2-40n-6}{6}v^3\right).
\end{split}\end{equation*}
The right-hand side is a polynomial in $v$ with nonnegative coefficients for every given $n\ge3$. Arguing similarly, one can easily check that, for $\Phi=B_n$, $n\ge 7$,
\begin{equation*}
\begin{split}
O_3(q)&=q^{n+3}+(2n-11)q^{n+2}+(n^2-12n+39)q^{n+1}+\dfrac{n^3-24n^2+185n-300}{6}q^n\\
&-\dfrac{n^3-15n^2+88n-176}{2}q^{n-1}+\dfrac{n^3-15n^2+76n-130}{2}q^{n-2}\\
&-\dfrac{n^3-21n^2+110n-174}{6}q^{n-3}-\dfrac{n^2-7n+12}{2}q^{n-4}\\
&=(v+1)^{n-4}\left(v^7+(2n-4)v^6+(n^2-6)v^5+\dfrac{n^3+6n^2+5n-60}{6}v^4\right.\\&+\left.\vphantom{\dfrac{n^3-37n+24}{6}v^4}\dfrac{n^3+9n^2-4n-42}{6}v^3+(n^2+n-5)v^2+(n-1)v\right).\\
\end{split}
\end{equation*}
}

In \cite{PakSoffer15} Pak and Soffer verified Higman's conjecture for $n\leq16$, and suggested that it probably fails for $n\geq59$. Similar conjectures are generalized to other finite groups with analogous structures, for example, finite pattern groups. Related problems are considered, and various tools from different aspects, including geometry, combinatorics, algebra, supercharacter theory have been developed. In this paper we prove that Isaacs' conjecture is true for certain generalizations of the Heisenberg group described in the next sections.

\newpage
\sect{Generalized Heisenberg group and its coadjoint orbits\label{sect:gener_Hei}}\addcontentsline{toc}{subsection}{\ref{sect:gener_Hei}. 
Generalized Heisenberg group and it's coadjoint orbits}
In this section we define a generalized Heisenberg group over a field and classify its coadjoint orbits when characteristic of the field is greater than 2 or equal to zero (Proposition~\ref{prop:gen_Hei_orbits}). In the case when the field is finite it is turned out that the Isaacs' conjecture is true (Corollary~\ref{coro:gen_Hei_Isaacs}).

Let $V$ be a vector space of dimension $n$ over a field $\Fp$, $\beta$ be an arbitrary bilinear form on $V$. We endow the set $V \times \Fp$ with a group structure according to the following rule: $$(v,\text{\,}r) \cdot (w,\text{\,}s) = (v + w,\text{\,}r + s + \beta(v, w)).$$ The identity and inverse elements are $(0,\text{\,}0)$ and $(-v,\text{\,}-\beta(v,v))$. We will call this group the \emph{generalized Heisenberg group} and denote it by $H_{\beta}$. 

It was shown in \cite{Deloup24}, that $H_{\widetilde{\beta}}$ and $H_{\beta}$ are equivalent as extensions of $V$ by $\Fp$ if and only if $\widetilde{\beta}(v, w) - \widetilde{\beta}(w, v) = \beta(v, w) - \beta(w, v)$ for all $v,w\in V$.

We choose an arbitrary basis in $V$. Then, as in \cite{Cushman24}, this group can be represented as a subgroup in the unitriangular matrices of size $(n+2)\times(n+2)$ as follows:
\[
(v,\text{\,} r) \mapsto 
\left(
\begin{array}{c@{\hspace{22pt}}c@{\hspace{15pt}}c}
1 & v & r \\
 & \ddots & \beta^{\wedge}(v) \\[5pt]
 &  & 1
\end{array}
\right),
\]
where other elements are zero, $\beta^{\wedge}(v) = \beta(v,\text{\,}\cdot\text{\,})$ and $v$ (respectively, $\beta^{\wedge}(v)$) is identified with the row (respectively, the column) of its coordinates in the chosen basis (respectively, in the dual basis). This can be easily verified by obvious matrix calculations. For evident reasons, it is a linear algebraic unipotent group defined over $\Fp$.

If $\text{char}(\Fp) \neq 2$, there is a logarithm mapping
\[\text{ln:}
\left(
\begin{array}{c@{\hspace{22pt}}c@{\hspace{15pt}}c}
1 & v & r \\
 & \ddots & \beta^{\wedge}(v) \\[5pt]
 &  & 1
\end{array}
\right)
\mapsto 
\left(
\begin{array}{c@{\hspace{25pt}}c@{\hspace{15pt}}c}
0 & v & r-\frac{1}{2}\beta(v,v) \\
 & \ddots & \beta^{\wedge}(v) \\[5pt]
 &  & 0
\end{array}
\right).
\]
It follows from the Baker--Campbell--Hausdorff formula that the logarithm image $\text{ln}(H_{\beta})$ is a Lie algebra with the commutator $[(v,\text{\,}r), (w,\text{\,}s)] = (0,\text{\,}-\beta(v,w)+\beta(w,v))$.  We will denote it by $\htt_{\beta}$.

The  matrices 
\[
e_i=\left(
\begin{array}{c@{\hspace{28pt}}c@{\hspace{15pt}}c}
0 & e_i & 0 \\
 & \ddots & \beta^{\wedge}(e_i) \\[5pt]
 &  & 0
\end{array}
\right), \text{  }1\leq i \leq n,\text{  } e_0=
\left(
\begin{array}{c@{\hspace{15pt}}c@{\hspace{15pt}}c}
0 & 0 & 1 \\
 & \ddots & 0 \\[5pt]
 &  & 0
\end{array}
\right)
\] form a basis in $\htt_{\beta}$. Using the trace form $\tr$, we identify $\htt^{*}_{\beta}$ with the space of all matrices of the form
\[
\left(
\begin{array}{c@{\hspace{15pt}}c@{\hspace{15pt}}c}
0 & 0 & 0 \\
v & \ddots & 0 \\[5pt]
r & 0 & 0
\end{array}
\right).
\]
\begin{flushleft}
$\text{Under this identification, the  vector $e^j$ from the dual basis in $\htt^{*}_{\beta}$ is }\left(
\begin{array}{c@{\hspace{15pt}}c@{\hspace{15pt}}c}
0 & 0 & 0 \\
e_j & \ddots & 0 \\[5pt]
0 & 0 & 0
\end{array}
\right) \text{for $j>0$  and}\break
\left(
\begin{array}{c@{\hspace{15pt}}c@{\hspace{15pt}}c}
0 & 0 & 0 \\
0 & \ddots & 0 \\[5pt]
1 & 0 & 0
\end{array}
\right) \text{  for $j = 0$.}$
\end{flushleft}

Let $(A)_C$ denote the matrix obtained from a matrix $A$ by replacing all elements except the second to $(n+1)$th elements of the first column with zero, i.e, 
\[
%\begin{pmatrix}
%a_{1,1} & a_{1,2} & \dots & a_{1,n+1} & %a_{1,n+2}\\
%a_{2,1} & a_{2,2} &  \dots & a_{2,n+1} & %a_{2,n+2}\\
%\vdots & \vdots & \ddots &\vdots & %\vdots\\
%a_{n+1,1} & a_{n+1,2} & \dots & %a_{n+1,n+1} & a_{n+1, n+2}\\
%a_{n+2,1} & a_{n+2,2} & \dots & %a_{n+2,n+1} & a_{n+2,n+2}
%\end{pmatrix}
(A)_C=
\begin{pmatrix}
0 & 0 & \dots & 0 & \quad0\quad\\
a_{2,1} & 0 &  \dots & 0 & \quad0\quad\\
\vdots & \vdots & \ddots & \vdots & \quad\vdots\quad\\
a_{n+1,1} & 0 & \dots & 0 & \quad0\quad\\
0 & 0 & \dots & 0 & \quad0\quad
\end{pmatrix}.
\]
\begin{flushleft}
$
\text{Similarly, }(A)_R=
\begin{pmatrix}
0 & 0\quad & \dots & 0 & 0\\
0 & 0\quad &  \dots & 0 & 0\\
\vdots & \vdots\quad & \ddots & \vdots & \vdots\\
0 & 0\quad & \dots & 0 & 0\\
0 & a_{n+2,2}\quad & \dots & a_{n+2,n+1} & 0
\end{pmatrix} \text{  and }(A)_{\boldsymbol{\bullet}}= 
\begin{pmatrix}
0 & 0 & \dots & 0 & \quad0\quad\\
0 & 0 &  \dots & 0 & \quad0\quad\\
\vdots & \vdots & \ddots & \vdots & \quad\vdots\quad\\
0 & 0 & \dots & 0 & \quad0\quad\\
a_{n+2,1} & 0 & \dots & 0 & \quad0\quad
\end{pmatrix}.
$
\end{flushleft}

By the definition of the coadjoint action, $(g.\lambda)(x) = \lambda(g^{-1}xg)$. So, first, $$\tr((g.\lambda)x) = \tr(((g\lambda g^{-1})_C + (g\lambda g^{-1})_{\boldsymbol{\bullet}} + (g\lambda g^{-1})_R)x)\text{,}$$ because
\begin{equation*}
\begin{split}
\tr((g\lambda g^{-1})_Cx + (g\lambda g^{-1})_{\boldsymbol{\bullet}}x + (g\lambda g^{-1})_Rx) &= \tr(g\lambda g^{-1}x)=\\ &= \tr(\lambda g^{-1}xg) = \lambda(g^{-1}xg) = (g.\lambda)(x) = \tr((g.\lambda)x).
\end{split}
\end{equation*}
Moreover, $$\tr((g\lambda g^{-1})_Cx + (g\lambda g^{-1})_{\boldsymbol{\bullet}}x + (g\lambda g^{-1})_Rx) = \langle(g\lambda g^{-1})_C,v\rangle + (g\lambda g^{-1})_{\boldsymbol{\bullet}}r + \langle(g\lambda g^{-1})_R,\beta^{\wedge}(v)\rangle,$$ where $\langle\text{\,}\cdot\text{\,},\text{\,}\cdot\text{\,}\rangle$ denotes the standard scalar product on $V$ with our basis, and $$\langle(g\lambda g^{-1})_R,\beta^{\wedge}(v)\rangle = \beta(v, (g\lambda g^{-1})_R) = \langle^{\wedge}\beta((g\lambda g^{-1})_R),v\rangle,$$ where $^{\wedge}\beta(f) = \beta(\text{\,}\cdot\text{\,},f)$. Thus, 
\begin{alignat*}{2}
\tr((g.\lambda)x) &= \\
&=\tr((g\lambda g^{-1})_Cx + (g\lambda g^{-1})_{\boldsymbol{\bullet}}x + (g\lambda g^{-1})_Rx) = \\
&=\langle(g\lambda g^{-1})_C,v\rangle + (g\lambda g^{-1})_{\boldsymbol{\bullet}}r + \langle^{\wedge}\beta((g\lambda g^{-1})_R),v\rangle =\\
&=\tr(((g\lambda g^{-1})_C + (g\lambda g^{-1})_{\boldsymbol{\bullet}} + {}^{\wedge}\beta((g\lambda g^{-1})_R))x).
\end{alignat*}
That is, $(g.\lambda)_C = (g\lambda g^{-1})_C + {}^{\wedge}\beta((g\lambda g^{-1})_R)$ and $(g.\lambda)_{\boldsymbol{\bullet}}=(g\lambda g^{-1})_{\boldsymbol{\bullet}}$.

Let the form $\beta$ has the Gram matrix $B$ in our basis in $V$. Then we can simply compute the explicit form of the $(g\lambda g^{-1})_C$, $(g\lambda g^{-1})_R$ and the $(g\lambda g^{-1})_{\boldsymbol{\bullet}}$:
\[
g\lambda g^{-1} =
\begin{pmatrix}
    \; 1 & v_1 & \dots & v_n & r\\
      & 1   &       &     & (B^{T}v)_1 \\
      &     & \ddots&     & \vdots \\
      &     &       & 1   & (B^{T}v)_n \\
      &     &       &     & 1 \\ 
\end{pmatrix}
\begin{pmatrix}
    0 &     &       &     &  \\
    w_1&  0   &       &     &  \\
    \vdots  &    &   \ddots  &     &  \\
     w_n &     &       &  0  &  \\
    s  &     &      &     &  0\;\\ 
\end{pmatrix}
\begin{pmatrix}
    1 & -v_1 & \dots & -v_n & v^TBv - r\\
      & 1   &       &     & -(B^{T}v)_1 \\
      &     & \ddots&     & \vdots \\
      &     &       & 1   & -(B^{T}v)_n \\
      &     &       &     & 1 \\ 
\end{pmatrix} =
\]
\[
=
\begin{pmatrix}
     *& * & \dots & *\quad\\
     w_1 + (B^{T}v)_1s & 0   &   \dots     & 0\quad \\
     \vdots &  \vdots   & \ddots &    \vdots \quad \\
     w_n + (B^{T}v)_ns &  0 &\dots& 0\quad \\
     s &   0  &    \dots      & 0\quad \\ 
\end{pmatrix}
\begin{pmatrix}
    1 & -v_1 & \dots & -v_n & v^TBv - r\\
      & 1   &       &     & -(B^{T}v)_1 \\
      &     & \ddots&     & \vdots \\
      &     &       & 1   & -(B^{T}v)_n \\
      &     &       &     & 1 \\ 
\end{pmatrix} =
\]
\[
=
\begin{pmatrix}
     *& * & \dots & * &*\quad\\
     w_1 + (B^{T}v)_1s & *   &   \dots     & * &*\quad \\
     \vdots &  \vdots   & \ddots &    \vdots & \vdots \quad\\
     w_n + (B^{T}v)_ns &  * &\dots& * &*\quad \\
     s &   -v_1s  &    \dots      & -v_ns &*\quad \\ 
\end{pmatrix}.
\]
So, according to the formula for $g.\lambda$ above, if $\lambda = (w, \text{\,}s)$,
\[
g.\lambda = (w + s(B^T-B)v, \text{\,}s).
\]

\newcommand{\slantquotient}[2]{\raisebox{1.5pt}{$#1$}/\raisebox{-1.5pt}{$#2$}}

Thus, we have proved
\prop{For\label{prop:gen_Hei_orbits} $\text{char}(\Fp) \neq 2$\textup, the classification of coadjoint orbits of the generalized Heisenberg group $H_{\beta}$ is as follows\textup:
\begin{equation*}
\begin{split}
&\text{\textup{(i)} if $s = 0$\textup, the orbit consists of a single point\textup;}\\
&\text{\textup{(ii)} if $s \neq 0$\textup, then the orbits characterized by vectors in ${V}/{\mathrm{Im}(B^T-B)}$.}
\end{split}
\end{equation*}}

In the case of finite field $\Fp = \Fp_{q}$, there are $q^n$ orbits consisting of a single point. There are $q^{\text{dim\,}\mathrm{Im}(B^T-B)}=q^{\text{rk}(B^T-B)}$ points in the orbits of the second type. The number of orbits of this dimension is $(q-1)q^{n-\text{rk}(B^T-B)}$. 

\[
\renewcommand{\arraystretch}{1.5}
\begin{tabular}{|c|c|c|}
    \hline
      &  The dimension &  The number of orbits \\
    \hline
     $s=0$ &  0 &  $q^n$ \\
    \hline
     $s \neq 0$ &  $\frac{1}{2}\text{dim\,}\mathrm{Im}(B^T-B) = \frac{1}{2}\text{rk}(B^T-B)$ &  $(q-1)q^{n-\text{rk}(B^T-B)}$ \\
    \hline
\end{tabular}
\]\\
This proves
\coro{The\label{coro:gen_Hei_Isaacs} Isaacs' conjecture holds for the generalized Heisenberg group $H_{\beta}$ over a finite field of characteristic greater than $2$.}
\exam{i) In the case where $\Fp = \Rp$ and $\beta$ is symplectic, the classification of orbits obtained in \cite{Cushman24} follows automatically from this result.

ii) In the case of the standard Heisenberg group
\[
\begin{pmatrix}
    \; 1 & x_1 & \dots & x_k & r\\
      & 1   &       &     & x_{k+1} \\
      &     & \ddots&     & \vdots \\
      &     &       & 1   & x_{2k} \\
      &     &       &     & 1 \\ 
\end{pmatrix}
\]
$V = \Fp^{2k}$ with coordinates $x_1, \text{\,} \dots\text{\,}, \text{\,} x_{2k}$ and by directly multiplying the matrices, one can verify that $\beta(x, y) = \displaystyle\sum_{i=1}^{k}x_iy_{i+k}$ and $\text{rk}(B^T-B)=2k$.
Thus, when $\Fp = \Fp_q$, the dimensions and numbers of orbits written out above coincide with well-known ones for the standard Heisenberg group.

iii) Also, in the case when $\beta$ is symplectic and $\Fp=\Fp_q$ is finite, one has $\text{rk}(B^T-B)=n$, and we immediately obtain the results of the papers \cite{DarafshehMisaghian08} and \cite{Misaghian11}.}

{\bf Remark.} According to \cite{Deloup24}, the standard Heisenberg group is isomorphic to the group $H_{\widetilde{\beta}}$ with symplectic $\widetilde{\beta}(x, y) = \dfrac{1}{2}\displaystyle\sum_{i=1}^{k}x_iy_{i+k} - \dfrac{1}{2}\displaystyle\sum_{i=k+1}^{2k}x_iy_{i-k}$, and thus case (ii) is a special case of case (iii).

\newpage

\sect{Classical Heisenberg-type groups}\label{Sfdcn}\addcontentsline{toc}{subsection}{\ref{Sfdcn}. Finite-dimensional case: notation}
In this section, we first briefly recall the definitions of classical finite-dimensional simple Lie algebras and fix the notation for the nilradicals of their Borel subalgebras and Heisenberg-type subgroups of these nilradicals.
%, which provide important class\break of Lie--Dynkin nil-algebras. 
This will be used to state the main results in latter sections. 
%we formulate our first main result, Theorem~\ref{mtheo:low_dim_fin}, which will be proved in the next section.

Pick $n\in\Zp_{>0}$. Let $\gt$ denote one of the Lie algebras $\slt_n(\Fp)$, $\sot_{2n}(\Fp)$, $\sot_{2n+1}(\Fp)$ or $\spt_{2n}(\Fp)$. We will assume till the end of the section that the characteristic of $\Fp$ is zero or greater than $\dim\gt$. The algebra $\sot_{2n}(\Fp)$ (respectively, $\sot_{2n+1}(\Fp)$ and $\spt_{2n}(\Fp)$) is realized as the sub\-al\-gebra of $\slt_{2n}(\Fp)$ (respectively, $\slt_{2n+1}(\Fp)$ and $\slt_{2n}(\Fp)$) consisting of all~$x$ such that $$\beta(u,xv)+\beta(xu,v)=0$$ for all $u,v$ in $\Fp^{2n}$ (respectively, in $\Fp^{2n+1}$ and $\Fp^{2n}$), where
\begin{equation*}
\beta(u,v)=\begin{cases}
\sum\nolimits_{i=1}^n(u_iv_{-i}+u_{-i}v_i)&\text{for }\sot_{2n}(\Fp),\\
u_0v_0+\sum\nolimits_{i=1}^n(u_iv_{-i}+u_{-i}v_i)&\text{for }\sot_{2n+1}(\Fp),\\
\sum\nolimits_{i=1}^n(u_iv_{-i}-u_{-i}v_i)&\text{for }\spt_{2n}(\Fp).
\end{cases}
\end{equation*}
Here for $\sot_{2n}(\Fp)$ (respectively, for $\sot_{2n+1}$ and $\spt_{2n}(\Fp))$ we denote by $e_1,\ldots,e_n,e_{-n},\ldots,e_{-1}$ (res\-pec\-tively, by $e_1,\ldots,e_n,e_0,e_{-n},\ldots,e_{-1}$ and $e_1,\ldots,e_n,e_{-n},\ldots,e_{-1}$) the standard basis of $\Fp^{2n}$ (res\-pec\-tively, of $\Fp^{2n+1}$ and $\Fp^{2n}$), and by $x_i$ the coordinate of a vector $x$ corresponding to $e_i$.

The set of all diagonal matrices from $\gt$ is a Cartan subalgebra of $\gt$; we denote it by $\htt$. Let $\Phi$ be the root system of $\gt$ with respect to $\htt$. Note that $\Phi$ is of type $A_{n-1}$ (respectively, $D_n$, $B_n$ and $C_n$) for $\slt_n(\Fp)$ (respectively, for $\sot_{2n}(\Fp)$, $\sot_{2n+1}(\Fp)$ and $\spt_{2n}(\Fp)$). The set of all upper-triangular matrices from $\gt$ is a Borel subalgebra of $\gt$ containing $\htt$; we denote it by $\bt$. Let $\Phi^+$ be the set of positive roots with respect to $\bt$. 
Denote by $\Pi$ the set of simple roots of $\Phi^+$. As usual \cite{Bou}, we identify $\Phi^+$ with the following subset of $\Rp^n$:
\begin{equation}\predisplaypenalty=0
\begin{split}
A_{n-1}^+&=\{\epsi_i-\epsi_j,~1\leq i<j\leq n\},\\
B_n^+&=\{\epsi_i-\epsi_j,~1\leq i<j\leq n\}\cup\{\epsi_i+\epsi_j,~1\leq i<j\leq n\}\cup\{\epsi_i,~1\leq i\leq n\},\\\label{formula:root_basis_nilradical_fin_dim}
C_n^+&=\{\epsi_i-\epsi_j,~1\leq i<j\leq n\}\cup\{\epsi_i+\epsi_j,~1\leq i<j\leq n\}\cup\{2\epsi_i,~1\leq i\leq n\},\\
D_n^+&=\{\epsi_i-\epsi_j,~1\leq i<j\leq n\}\cup\{\epsi_i+\epsi_j,~1\leq i<j\leq n\}.\\
\end{split}
\end{equation}
Here $\{\epsi_i\}_{i=1}^n$ is the standard basis of $\Rp^n$.

Denote by $\nt$ the algebra of all strictly upper-triangular matrices from $\gt$. Then $\nt$ has a basis consisting of root vectors $e_{\alpha}$, $\alpha\in\Phi^+$, where
\begin{equation*}\predisplaypenalty=0
\begin{split}
e_{\epsi_i}&=\sqrt{2}(e_{0,i}-e_{-i,0}),~e_{2\epsi_i}=e_{i,-i},\\
e_{\epsi_i-\epsi_j}&=\begin{cases}
e_{i,j}&\text{for }A_{n-1},\\
e_{i,j}-e_{-j,-i}&\text{for }B_n,~C_n\text{ and }D_n,
\end{cases}\\
e_{\epsi_i+\epsi_j}&=\begin{cases}
e_{i,-j}-e_{j,-i}&\text{for }B_n\text{ and }D_n,\\
e_{i,-j}+e_{j,-i}&\text{for }C_n,
\end{cases}
\end{split}
\end{equation*}
and $e_{i,j}$ are the usual elementary matrices. For $\sot_{2n}(\Fp)$ (respectively, for $\sot_{2n+1}(\Fp)$ and $\spt_{2n}(\Fp)$) we index the rows (from left to right) and the columns (from top to bottom) of matrices by the numbers $1,\ldots,n,-n,\ldots,-1$ (respectively, by the numbers $1,\ldots,n,0,-n,\ldots,-1$ and $1,\ldots,n,-n,\ldots,-1$). Note that $$\gt=\htt\oplus\nt\oplus\nt_-,$$ where $\nt_-=\langle e_{-\alpha},~\alpha\in\Phi^+\rangle_{\Fp}$, and, by definition, $e_{-\alpha}=e_{\alpha}^T$. (The superscript~$T$ always indicates matrix transposition.) The set $\{e_{\alpha},~\alpha\in\Phi\}$ can be uniquely extended to a Chevalley basis of $\gt$. Clearly, $\nt$ is the nilradical of the Borel subalgebra $\bt$ and $e_\alpha, \alpha\in\Phi^+,$ is a basis of $\nt^*$.

We denote $N=\exp(\nt)$, so that $\nt$ is the Lie algebra of the algebraic group $N$. The group $N$ acts on $\nt$ via the adjoint action, and the dual action of $N$ on $\nt^*$ is called \emph{coadjoint}. 
We denote by $e_\alpha^*$, $\alpha\in\Phi^+,$ the standard dual basis to $\nt^*$. 

\begin{comment}
\defi{Roots $\alpha, \gamma\in\Phi^+$ are called
$\beta$-\emph{singular} if $\alpha+\gamma=\beta$. The set of all $\beta$-singular roots is denoted by $\Sing(\beta)$.}
It is easy to see that
\begin{equation}\label{Esroot}
\begin{split}
\Sing(\epsi_i-\epsi_j)=&\bigcup_{l=i+1}^{j-1}\{\epsi_i-\epsi_l,\epsi_l-\epsi_j\},~1\leq i < j\leq n,\\
\Sing(\epsi_i)=&\bigcup_{l=i+1}^{n}\{\epsi_i-\epsi_l,\epsi_l\},~1\leq i\leq n,\\
\Sing(\epsi_i+\epsi_j)=&\bigcup_{l=i+1}^{j-1}\{\epsi_i-\epsi_l,\epsi_l+\epsi_j\}
\cup\bigcup_{l=j+1}^n\{\epsi_i-\epsi_l,\epsi_j+\epsi_l\}\cup\\
&\bigcup_{l=j+1}^n\{\epsi_i+\epsi_l,\epsi_j-\epsi_l\}\cup S_{ij}
,~1\leq i\le j\leq
n,\mbox{ where}\\
S_{ij}=&\begin{cases}\{\epsi_i, \epsi_j\},&\text{if }\Phi=B_n,\\
\{\epsi_i-\epsi_j, 2\epsi_j\},&\text{if }i\ne j,~\Phi=C_n,\\
\varnothing,&\text{if }\Phi=D_n
\end{cases}.
\end{split}
\end{equation}
\end{comment}

Recall that there is a natural partial order on $\Phi^+$ defined as follows: 
for $\alpha, \beta\in\Phi$ we set $\alpha\le\beta$ (or $\beta\ge\alpha$) iff $\beta-\alpha$ is a sum of vectors from $\Phi^+$ or if $\alpha=\beta$. 
%$\alpha\leq\beta$ if $\beta-\alpha$ is a sum of positive roots.
For example, $\epsi_i-\epsi_j\leq\epsi_a-\epsi_b$ if and only if $a\leq i$ and $j\leq b$; one can easily obtain similar conditions for other types of roots and root systems.

{\bf Remark.} It is convenient to draw schematically root systems under consideration as follows. 
We will draw $\Phi=A_{n-1}$ as the lower-triangular chessboard $n\times n$, where a root $\alpha=\epsi_i-\epsi_j$, $i<j$, corresponds to the box $(j,i)$. The root system $\Phi=C_n$ is drawn as a half of lower-triangular chessboard $2n\times2n$. 
The choice ``lower-triangular" for the upper-triangular matrices and the chessboard is used because we are mainly interested in the coadjoint representation (which is somehow symmetric to the adjoint representation). 

Precisely, we will index the rows and the columns of the chessboard by the indices 
$$1,\ldots,n,-n,\ldots,-1.$$
A root $\alpha=\epsi_i-\epsi_j$, $i<j$, will be drawn as the box $(j,i)$; a root $\alpha=\epsi_i+\epsi_j$ corresponds to the box $(j,-i)$. (Similar convention is applied for $\Phi=D_n$.) A root $\alpha=2\epsi_i$ corresponds to the box $(-i,i)$. Finally, $\Phi=B_n$ will be drawn as a half of the lower-triangular chessboard $(2n+1)\times(2n+1)$. We will enumerate rows and columns by the indices 
$$1,\ldots,n,0,-n,\ldots,-1.$$
Roots of the form $\alpha=\epsi_i\pm\epsi_j$ satisfy the same conditions as above, while a root $\alpha=\epsi_i$ corresponds to the box $(0,i)$.

The corresponding pictures for $A_4, B_4, C_4, D_4$ are given below.
$$\begin{array}{c}A_4:
~{\Autonumfalse\mymatrix{
\lNote{1}\Note{1}\Bot{2pt}& \Note{2}\pho&\Note{3}\pho&\Note{4}\pho&\\
\lNote{2}\gray\pho\Rt{2pt}&\Note{2}\Bot{2pt}\pho&\pho&&\\
\lNote{3}\gray\pho&\gray \Rt{2pt}&\Bot{2pt}&&\\
\lNote{4}\gray&\gray\star&\gray \Rt{2pt}&\Bot{2pt}&\\
\lNote{5}\gray&\gray&\gray &\gray\Rt{2pt}&\pho\\
}}\\
\star\to\epsi_2-\epsi_4
\end{array},
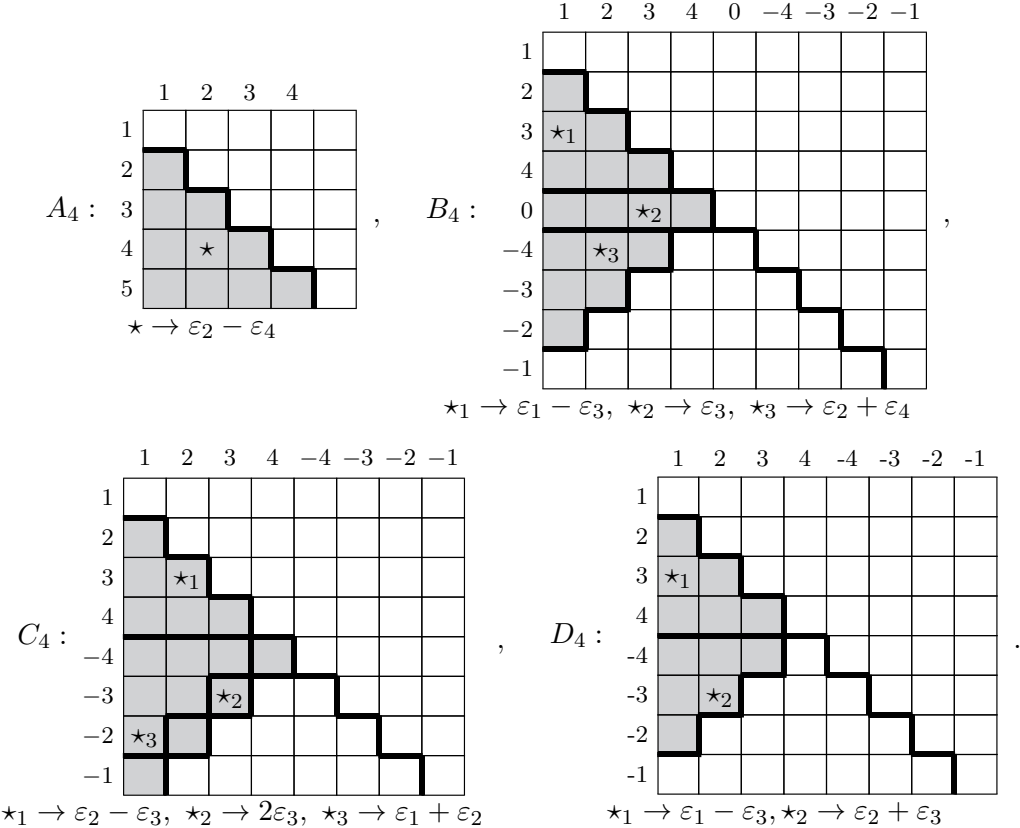
\hspace{10pt}
\begin{array}{c}B_4:
~{\Autonumfalse\mymatrix{
\lNote{1}\Note{1}\Bot{2pt}\pho& \Note{2}\pho& \Note{3}\pho& \Note{4}\pho& \Note{0}\pho& \Note{$-4$}\pho& \Note{$-3$}\pho& \Note{$-2$}\pho& \Note{$-1$}\pho\\
\lNote{2}\gray\Rt{2pt}\pho& \Bot{2pt}\pho& \pho& \pho& \pho& \pho& \pho& \pho& \pho\\
\lNote{3}\gray\star_1& \Rt{2pt}\gray\pho& \Bot{2pt}\pho& \pho& \pho& \pho& \pho& \pho& \pho\\
\lNote{4}\gray\Bot{2pt}\pho& \gray\Bot{2pt}\pho& \gray\Rt{2pt}\Bot{2pt}\pho& \Bot{2pt}\pho& \pho& \pho& \pho& \pho& \pho\\
\lNote{0}\gray\Bot{2pt}\pho& \gray\Bot{2pt}\pho& \gray\Bot{2pt}\star_2& \gray\Rt{2pt}\Bot{2pt}\pho& \Bot{2pt}\pho& \pho& \pho& \pho& \pho\\
\lNote{$-4$}\gray\pho& \gray\star_3& \Bot{2pt}\Rt{2pt}\gray\pho& \pho& \Rt{2pt}\pho& \Bot{2pt}\pho& \pho& \pho& \pho\\
\lNote{$-3$}\gray\pho& \gray\Bot{2pt}\Rt{2pt}\pho& \pho& \pho& \pho& \Rt{2pt}\pho& \Bot{2pt}\pho& \pho& \pho\\
\lNote{$-2$}\gray\Bot{2pt}\Rt{2pt}\pho& \pho& \pho& \pho& \pho& \pho& \Rt{2pt}\pho& \Bot{2pt}\pho& \pho\\
\lNote{$-1$}\pho& \pho& \pho& \pho& \pho& \pho& \pho& \Rt{2pt}\pho& \pho\\}}\\
\star_1\to\epsi_1-\epsi_3,~\star_2\to\epsi_3,~\star_3\to\epsi_2+\epsi_4
\end{array},\hspace{10pt}$$
$$\begin{array}{c}C_4:
{\Autonumfalse\mymatrix{
\lNote{1}\Note{1}\Bot{2pt}\pho& \Note{2}\pho& \Note{3}\pho& \Note{4}\pho& \Note{$-4$}\pho& \Note{$-3$}\pho& \Note{$-2$}\pho& \Note{$-1$}\pho\\
\lNote{2}\gray\Rt{2pt}\pho& \Bot{2pt}\pho& \pho& \pho& \pho& \pho& \pho& \pho\\
\lNote{3}\gray\pho& \gray\Rt{2pt}\star_1& \Bot{2pt}\pho& \pho& \pho& \pho& \pho& \pho\\
\lNote{4}\gray\Bot{2pt}\pho& \gray\Bot{2pt}& \gray\Bot{2pt}\Rt{2pt}\pho& \Bot{2pt}\pho& \pho& \pho& \pho& \pho\\
\lNote{$-4$}\gray\pho& \gray\pho& \gray\Bot{2pt}\Rt{2pt}& \gray\Rt{2pt}& \Bot{2pt}\pho& \pho& \pho& \pho\\
\lNote{$-3$}\gray\pho& \gray\Bot{2pt}\Rt{2pt}\pho& \gray\star_2& \Lft{2pt}\Top{2pt}\pho& \Rt{2pt}\pho& \Bot{2pt}\pho& \pho& \pho\\
\lNote{$-2$}\gray\Bot{2pt}\Rt{2pt}\star_3& \gray\pho& \Lft{2pt}\Top{2pt}\pho& \pho& \pho& \Rt{2pt}\pho& \Bot{2pt}\pho& \pho\\
\lNote{$-1$}\gray\pho& \Lft{2pt}\Top{2pt}\pho& \pho& \pho& \pho& \pho& \Rt{2pt}\pho& \pho\\}}\\
\star_1\to\epsi_2-\epsi_3,~\star_2\to2\epsi_3,~\star_3\to\epsi_1+\epsi_2
\end{array},\hspace{10pt}
\begin{array}{c}D_4:
~{\Autonumfalse
\mymatrix{
\lNote{1} \Note{1}\Bot{2pt}& \Note{2}\pho& \Note{3}\pho& \Note{4}\pho& \Note{-4}\pho& \Note{-3}\pho& \Note{-2}\pho& \Note{-1}\pho\\
\lNote{2} \gray\Note{1}\Rt{2pt}& \Note{2}\Bot{2pt}\pho& \Note{3}\pho&&&&&\\
\lNote{3}\gray\star_1&\gray \Rt{2pt}&\Bot{2pt}\pho&&&&&\\
\lNote{4}\gray\Bot{2pt}&\gray\Bot{2pt}&\gray\Bot{2pt}\Rt{2pt}&\Bot{2pt}&&&&\\
\lNote{-4}\gray&\gray&\gray\Rt{2pt}\Bot{2pt}&\Rt{2pt}&\Bot{2pt}&&&\\
\lNote{-3}\gray&\gray\star_2 \Rt{2pt}\Bot{2pt}&&&\Rt{2pt}&\Bot{2pt}&&\\
\lNote{-2}\gray\Bot{2pt}\Rt{2pt}& \pho& \pho&&&\Rt{2pt}&\Bot{2pt}&\\
\lNote{-1}& \pho& \pho&&&&\Rt{2pt}&\\
}}\\
\star_1\to\epsi_1-\epsi_3, \star_2\to\epsi_2+\epsi_3\end{array}.$$
In these pictures the cells corresponding to the roots of $\Phi^+$ are marked by gray color and for several such cells we provide the explicit presentation of the corresponding roots

We denote
$$  
\begin{array}{ll}
&\col\colon\Phi^+\to\{1,\ldots,n\}\colon\col(\epsi_i\pm\epsi_j)=\col(\epsi_i)=i,\\
&\row\colon\Phi^+\to\{-n,\ldots,n\}\colon\row(\epsi_i\pm\epsi_j)=\mp j,\row(\epsi_i)=0,\\
\end{array}
$$
and for arbitrary $-n+1\leq i\leq n-1$ and $1\leq j\leq n$ the sets
$$
\begin{array}{ll}
&\Ro_i = R_i(\Phi) = \{\alpha\in\Phi^+\mid \row(\alpha)=i\},\\
&\Co_j = C_j(\Phi) = \{\alpha\in\Phi^+\mid \col(\alpha)=j\}\\
\end{array}
$$
are called the \emph{$i$th row} and the \emph{$j$th column} of $\Phi^+$ respectively.

The standard Heisenberg group can be defined as follows. Let $\Phi=A_{n-1}$, so that $\nt$ is the Lie algebra of strictly upper-triangular matrices with zeros on the diagonal, and $N=U_n(\Fp)$ is the unitriangular group, i.e., the group of upper-triangular matrices with units on the diagonal. Then one can define the Heisenberg Lie algebra $$\hei_n=\langle e_{\alpha},~\alpha\in\Ro_1\cup\Co_n\rangle_{\Fp}$$
and the Heisenberg group $$\Hei_n=\exp\hei_n=\{g\in N\mid g_{i,j}=0\text{ if }i<j\text{ and }\epsi_i-\epsi_j\notin\Ro_1\cup\Co_n\}.$$
Next, one can define so-called \emph{$2$-layered Heisenberg group} \cite{Nien20} as $$\Hei_{n,2}=\exp\hei_{n,2},$$ where $$\hei_{n,2}=\langle e_{\alpha},~\alpha\in\Ro_1\cup\Ro_2\cup\Co_n\cup\Co_{n-1}\rangle_\Fp.
$$
This group is a natural generalization of the usual Heisenberg group. It is a so-called pattern group, which makes it interesting from several points of view. For example, in \cite{Nien20}, the number of irreducible representations of $\Hei_{n,2}$ was computed.

Similarly, for each classical root system $\Phi$ of type $B_n$, $C_n$ or $D_n$, one can define the corresponding 2-layered Heisenberg-type group as follows: $$\Hei_2(\Phi)=\exp\hei_2(\Phi),$$ where
$$\hei_2(\Phi)=\langle e_{\alpha},~\alpha\in\Co_1\cup\Co_2\rangle_{\Fp}.$$
(Since the matrices from $\nt$ are either symmetric or antisymmetric, it is natural to use just columns without rows.) For example, we schematically drew the algebra $\hei_2(\Phi)$ for $\Phi=C_4$ and $D_4$ (boxes corresponding to the roots from $\Co_1\cup\Co_2$ are marked by gray):

$$\begin{array}{c}C_4:
{\Autonumfalse\mymatrix{
\lNote{1}\Note{1}\Bot{2pt}\pho& \Note{2}\pho& \Note{3}\pho& \Note{4}\pho& \Note{$-4$}\pho& \Note{$-3$}\pho& \Note{$-2$}\pho& \Note{$-1$}\pho\\
\lNote{2}\gray\Rt{2pt}\pho& \Bot{2pt}\pho& \pho& \pho& \pho& \pho& \pho& \pho\\
\lNote{3}\gray\pho& \gray\Rt{2pt}\pho& \Bot{2pt}\pho& \pho& \pho& \pho& \pho& \pho\\
\lNote{4}\gray\Bot{2pt}\pho& \gray\Bot{2pt}&\Bot{2pt}\Rt{2pt}\pho& \Bot{2pt}\pho& \pho& \pho& \pho& \pho\\
\lNote{$-4$}\gray\pho& \gray\pho& \Bot{2pt}\Rt{2pt}& \Rt{2pt}& \Bot{2pt}\pho& \pho& \pho& \pho\\
\lNote{$-3$}\gray\pho& \gray\Bot{2pt}\Rt{2pt}\pho& \pho& \Lft{2pt}\Top{2pt}\pho& \Rt{2pt}\pho& \Bot{2pt}\pho& \pho& \pho\\
\lNote{$-2$}\gray\Bot{2pt}\Rt{2pt}\pho& \gray\pho& \Lft{2pt}\Top{2pt}\pho& \pho& \pho& \Rt{2pt}\pho& \Bot{2pt}\pho& \pho\\
\lNote{$-1$}\gray\pho& \Lft{2pt}\Top{2pt}\pho& \pho& \pho& \pho& \pho& \Rt{2pt}\pho& \pho\\}}
\end{array},
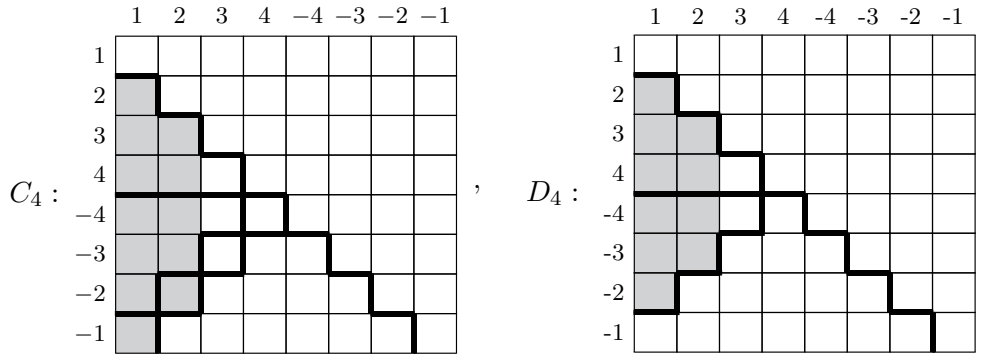
\hspace{10pt}
\begin{array}{c}D_4:
~{\Autonumfalse
\mymatrix{
\lNote{1} \Note{1}\Bot{2pt}& \Note{2}\pho& \Note{3}\pho& \Note{4}\pho& \Note{-4}\pho& \Note{-3}\pho& \Note{-2}\pho& \Note{-1}\pho\\
\lNote{2} \gray\Note{1}\Rt{2pt}& \Note{2}\Bot{2pt}\pho& \Note{3}\pho&&&&&\\
\lNote{3}\gray\pho&\gray \Rt{2pt}&\Bot{2pt}\pho&&&&&\\
\lNote{4}\gray\Bot{2pt}&\gray\Bot{2pt}&\Bot{2pt}\Rt{2pt}&\Bot{2pt}&&&&\\
\lNote{-4}\gray&\gray&\Rt{2pt}\Bot{2pt}&\Rt{2pt}&\Bot{2pt}&&&\\
\lNote{-3}\gray&\gray\pho \Rt{2pt}\Bot{2pt}&&&\Rt{2pt}&\Bot{2pt}&&\\
\lNote{-2}\gray\Bot{2pt}\Rt{2pt}& \pho& \pho&&&\Rt{2pt}&\Bot{2pt}&\\
\lNote{-1}& \pho& \pho&&&&\Rt{2pt}&\\
}}\end{array}.$$

In the next section, we calculate the number of irreducible representations of $\Hei_2(\Phi)$ for $\Phi=B_n$ or $D_n$ and prove that the Isaacs conjecture holds in these cases, see Theorem~\ref{theo:2_Hei_Isaacs}, which we consider as the second main result of the paper.

\newpage

\sect{The Mackey method for orthogonal 2-layered Heisenberg groups}\label{sect:Mackey}\addcontentsline{toc}{subsection}{\ref{sect:Mackey}. The Mackey method for classical 2-layered Heisenberg groups over finite field}

Let $\Gr$ be an arbitrary finite group, $A$, $B$ be its subgroups such that $\Gr$ is the semi-direct product of $A$ and $B$, i.e., $\Gr=AB$, $A\triangleleft\Gr$ and
$A\cap B=\{1\}$ (we denote $\Gr=A\rtimes B$). Assume additionally that $A$ is an abelian group and that $\psi$ is its irreducible character. The \emph{centralizer}
of $\psi$ in the group $B$ (or the \emph{little group}) is the subgroup of the form $$B^{\psi}=\{b\in
B\mid\psi\circ\tau_b=\psi\},$$ where $\tau_b\colon A\to A.~a\mapsto b^{-1}ab$. Clearly, $B$ acts on the set of irreducible characters of $A$ by the formula $b.\psi=\psi\circ\tau_b$, and $B^{\psi}$ is the stabilizer of $\psi$ under this action.

An arbitrary element $g\in\Gr$ can be uniquely represented as
$g=ab$, where $a\in A, b\in B$; this defines the maps
$$\pi^{\Gr}_{A}\colon\Gr\to A\colon g\mapsto a,~\pi^{\Gr}_{B}\colon\Gr\to B\colon g\mapsto b,~g\in\Gr.$$ Note that $\pi^{\Gr}_B$ is a group homomorphism, while $\pi^{\Gr}_A$, in general, is not. For an arbitrary subgroup $\wt B$ in~$B$ and arbitrary characters
$\psi\colon A\to\Cp$ and $\eta\colon\wt B\to\Cp$, one can define the character $\psi_0\eta_0$ of the group
$A\rtimes\wt B=A\wt B$, where
$$
\psi_0=\psi\circ\pi^{A\rtimes\wt
B}_{A},~\eta_0=\eta\circ\pi^{A\rtimes\wt B}_{\wt B}.
$$
Denote by $\Irr{\Gr}$ the set of irreducible characters of $\Gr$. The basic idea of the Mackey method, which allows us
to reduce the study of representations of the group $\Gr$ to the study of representations of
its subgroups $A$ and $B$, can be formulated as follows. \theo{Let $\Gr=A\rtimes B$
be a finite group with $A$ being an abelian group. Then each irreducible character $\chi$ of the group $\Gr$ has the form
\begin{equation}
\chi=\chi_{\psi,\eta}=\Ind{A\rtimes B^{\psi}}{G}{\psi_0\eta_0}
\label{formula_semi_direct}
\end{equation}
for certain $\psi\in\Irr{A}$\textup, $\eta\in\Irr{B}$. Conversely\textup, every function of the form \textup{(\ref{formula_semi_direct})} is an irreducible character of $\Gr$. Furthermore\textup, let $\{\psi_i\}$ be the set of representatives of $B$-orbits on $\Irr{A}$\textup, then\break $\chi_{\psi_i,\eta}=\chi_{\psi_j,\eta'}$ if and only if $i=j$ and $\eta=\eta'$.}{\cite[8.2, Proposition
25]{Serre77}}\label{theo_semi_direct}

For the root system $B_{n+1}$ (respectively, the root system $D_{n+1}$), let $\at_1$ (respectively, $\at_2$) be the space $\langle e_{\alpha},~\alpha\in\Co_1\rangle_{\Fp}$ and $\bt_1$ (respectively, $b_2$) be the space $\langle e_{\alpha},~\alpha\in\Co_2\rangle_{\Fp}$. It is obvious that $\at_1 \oplus \bt_1 = \sot_{2n+3}(\Fp)$, $\at_2 \oplus \bt_2 = \sot_{2n+2}(\Fp)$ as vector spaces and $\at_1 \triangleleft \sot_{2n+3}(\Fp)$, $\at_2 \triangleleft \sot_{2n+2}(\Fp)$ are abelian ideals. From now on, we will assume that the characteristic of the field is greater than 3 or is zero. So, $\Hei_{2}(B_{n+1}) = A_1 \rtimes B_1$ and $\Hei_{2}(D_{n+1}) = A_2 \rtimes B_2$, where $A_i = \exp\at_i$ and $B_i = \exp\bt_i$, $i=1,2$. Thus, over a finite field the conditions of the Mackey method are met for types $B$ and $D$. For type $C$, the Mackey method is not directly applicable, because the ideal $\langle e_{\alpha},~\alpha\in\Co_1\rangle_{\Fp}$ is not abelian. From now on, we will omit the indices $1$ and $2$ where it does not matter.

It is elementary to verify that the groups $A_1$ and $A_2$ consist of matrices of the form
\[
\begin{pmatrix}
    1 & a_1 & \dots & a_n & a_0 & a_{-n} & \dots & a_{-1} & -\displaystyle\sum a_ia_{-i} \\
   & & & & & & & & -a_{-1}\\
   & &  & & & & & & \vdots \\
   & &  & &  & & & & -a_{-n}\\
   & &  & &  & & & & -a_0\\
   & &  & & &  & & & -a_n \\
   & &  & & & & & & \vdots\\
   & &  & & & & &  & -a_1\\
   & &  & & & & & & 1
\end{pmatrix} \text{ and}
\]
\[
\begin{pmatrix}
    1 & a_1 & \dots & a_n & a_{-n} & \dots & a_{-1} & -\displaystyle\sum a_ia_{-i} \\
    & & & & & & & -a_{-1}\\
    &  & & & & & & \vdots \\
    &  & &  & & & & -a_{-n}\\
    &  & & &  & & & -a_n \\
    &  & & & & & & \vdots\\
    &  & & & & &  & -a_1\\
    &  & & & & & & 1
\end{pmatrix}
\]
respectively, with units on the diagonal and zeros in other places, and the groups $B_1$ and $B_2$ consists of matrices of the form
\[
\begin{pmatrix}
    1& 0 & 0 & \dots & 0 & 0 & 0 & \dots & 0 & 0 & 0 \\
     & 1 & b_2 & \dots & b_{n} & b_0&  b_{-n} & \dots & b_{-2} & -\displaystyle\sum b_ib_{-i} & 0 \\
    && & & & & & & & -b_{-2} & 0\\
    &&  & & & & & & & \vdots & \vdots\\
    &&  & &  & & & & & -b_{-n} & 0\\
    &&  & &  & & & & & -b_0 & 0\\
    &&  & & &  & & & & -b_{n} & 0\\
    &&  & & & & & & & \vdots & \vdots\\
    &&  & & & & &  & & -b_2 & 0\\
    &&  & & & & & & & 1 & 0 \\
    &&  & & & & & & & & 1 \\
\end{pmatrix}\text{ and}
\]
\[
\begin{pmatrix}
    1& 0 & 0 & \dots & 0 & 0 & \dots & 0 & 0 & 0 \\
     & 1 & b_2 & \dots & b_{n} & b_{-n} & \dots & b_{-2} & -\displaystyle\sum b_ib_{-i} & 0 \\
    & & & & & & & & -b_{-2} & 0\\
    &  & & & & & & & \vdots & \vdots\\
    &  & &  & & & & & -b_{-n} & 0\\
    &  & & &  & & & & -b_{n} & 0\\
    &  & & & & & & & \vdots & \vdots\\
    &  & & & & &  & & -b_2 & 0\\
    &  & & & & & & & 1 & 0 \\
    &  & & & & & & & & 1 \\
\end{pmatrix},
\]
respectively, where, unless otherwise specified, we assume by default that the summation runs over all admissible indices $i>0$.

It is easy to see that the group $A_i$ (respectively, the group $B_i$) is isomorphic to $(\Fp^{k_i}, +)$ (respectively, to $(\Fp^{k_i-2}, +)$), $i=1,2$, $k_1=2n+1$, $k_2 = 2n$, where we associate to an element $a$ (respectively, $b$) the vector with coordinates $(a_1, \ldots, a_{-1}) \in \Fp^{k}$ (respectively, $(b_2, \ldots , b_{-2}) \in \Fp^{k-2}$). It is obvious that $b^{-1}$ can be obtained from $b$ by replacing $b_j$ by $-b_j$ for all $j$. We can simply compute $bab^{-1}$ by multiplying three matrices. Since the groups $A$ are normal and, by construction of the isomorphisms $A$ to $\Fp^{k}$, we are only interested in the second to the $(2n+2)$th for $B_{n+1}$ (and in the second to the $(2n+1)$th for $D_{n+1}$) elements of the first row of the matrix $bab^{-1}$. The first row of $bab^{-1}$ has the form $$(1,a_1, -a_1b_2 + a_2, \dots , -a_1b_k+a_k, \dots, -a_1b_{-2}+a_{-2}, -a_1\displaystyle\sum b_ib_{-i} + \displaystyle\sum_{i \neq \pm 1}a_ib_{-i}+a_{-1}, *),$$ where $*$ denotes some element of the field $\Fp$. Thus, identifying $A$ with $\Fp^k$, we have 
$$bab^{-1} = (a_1, -a_1b_2 + a_2, \dots , -a_1b_k+a_k, \dots, -a_1b_{-2}+a_{-2}, -a_1\displaystyle\sum b_ib_{-i} + \displaystyle\sum_{i \neq \pm 1}a_ib_{-i}+a_{-1}).$$ 

We will assume till the end of the section that the $\Fp = \Fp_q$ is finite and $\chara\Fp_q = p > 3$. The irreducible characters of $(\Fp^{k}_q, +)$ are well known:
$$\forall a \in \Fp^{k}_q~~\psi_x(a) = e^{\frac{2\pi i}{p}\text{Tr}_{\Fp_q/\Fp_p}(\langle x, a\rangle)},$$ where $\text{Tr}_{\Fp_q/\Fp_p}(\alpha) = \alpha + \alpha^p + \alpha^{p^2} + \cdots + \alpha^{p^{n-1}}$ and $\langle\text{\,}\cdot\text{\,},\text{\,}\cdot\text{\,}\rangle$ denotes the standard bilinear form on $\Fp^{k}_q$: $$\langle x,y\rangle=x_1y_1+\ldots+x_ky_k,~x,y\in\Fp_q^k,$$ see, for example, \cite[Theorem 5.7 and Section 9.1]{LidlNiederreiter09} for more details. 

By definition, $b \in B^{\psi}$ if and only if $\psi_x(bab^{-1}) = \psi_x(a)$ for all $a \in \Fp^{k}_q$, hence, $$\text{Tr}_{\Fp_q/\Fp_p}(\langle x, a\rangle) = \text{Tr}_{\Fp_q/\Fp_p}(\langle x, bab^{-1}\rangle).$$ The map $\mathrm{Tr}$ is $\Fp_p$-linear \cite{LidlNiederreiter09}, so we have $$\text{Tr}_{\Fp_q/\Fp_p}(\langle x, a - bab^{-1}\rangle) = 0\text{ for all }a \in \Fp^{k}_q.$$ Also, $$ a - bab^{-1} = (0, a_1b_2, \ldots , a_1b_k, \ldots, a_1b_{-2}, a_1\displaystyle\sum b_ib_{-i} - \displaystyle\sum_{i \neq \pm 1}a_ib_{-i}).$$
We will consider three different cases of $x$.

i) If $x = (x_1, 0, \dots, 0) \in \Fp^{k}_q,~~x_1\in\Fp_q$, then $B^{\psi} = B$. By the formula for characters of $(\Fp^{k}, +)$ above, there are $q^{k-2}$ different irreducible characters of the group $B$.

ii) If $x = (x_1, \dots, x_{-2}, 0)$ and the set $J \overset{\mathrm{def}}{=} \{l\mid l\neq\pm1~\text{and}~x_l \neq 0\}$ is non-empty, we have $$\text{Tr}_{\Fp_q/\Fp_p}(\displaystyle\sum_{j \in J}a_1x_jb_j)=\text{Tr}_{\Fp_q/\Fp_p}(a_1(\displaystyle\sum_{j \in J}x_jb_j))=0\text{ for all } a_1 \in \Fp_q.$$ By \cite[Theorem 2.24]{LidlNiederreiter09}, it is equal to $\displaystyle\sum_{j \in J}x_jb_j = 0.$ That is the group $B^{\psi}$ is the vector subspace of codimension $1$ in $\Fp^{k-2}_q$ defined by the linear equation $\displaystyle\sum_{i \neq \pm 1}x_ib_i = 0$. By the same reasons as for $B$ above, the group $B^{\psi}$ has $q^{k-3}$ irreducible characters and there are $\dfrac{q^{k-2} - 1}{q - 1}$ different groups $B^{\psi}$.

iii) If $x_{-1} \neq 0$, then, putting $a_1=0$ we have $$\text{Tr}_{\Fp_q/\Fp_p}(\displaystyle\sum_{i \neq \pm 1}a_ix_{-1}b_{-i})=0\text{ for all }a_2, \dots, a_{-2} \in \Fp_q.$$ By setting all $a_i$, except $a_l$, to zero, for all $l$, $l \neq \pm 1$, we obtain by \cite[Theorem 2.24]{LidlNiederreiter09} that $b_i = 0$ for all $i$, and $B^{\psi} = \{e\}$ has only one character.

Let us find the orbit of an arbitrary irreducible character $\psi_x$ under the action of the group $B$. One has $b^{-1}.\psi_x = \psi_y$ if and only if $\psi_x(bab^{-1}) = \psi_y(a)$ for all $a \in A$, which is equivalent to the condition $\text{Tr}_{\Fp_q/\Fp_p}(\langle x, bab^{-1}\rangle) = \text{Tr}_{\Fp_q/\Fp_p}(\langle y, a\rangle)$ for all $a \in \Fp^{k}_q$. Putting $a_1=0$, we have $$bab^{-1} = (0, a_2, \ldots , a_k, \ldots, a_{-2}, \displaystyle\sum_{i \neq \pm 1}a_ib_{-i}+a_{-1})~~\text{and}$$ 
\begin{equation*} 
\begin{split}
    \text{Tr}_{\Fp_q/\Fp_p}(\langle x, bab^{-1}\rangle) &= \text{Tr}_{\Fp_q/\Fp_p}(\displaystyle\sum_{i \neq \pm 1}x_ia_i + x_{-1}\displaystyle\sum_{i \neq \pm 1}a_ib_{-i}+x_{-1}a_{-1}) = \\
    & =\text{Tr}_{\Fp_q/\Fp_p}(\displaystyle\sum_{i \neq \pm 1}y_ia_i+y_{-1}a_{-1})\text{ for all }a_2, \dots, a_{-1} \in \Fp_q. 
\end{split} 
\end{equation*} 
Since the map Tr is linear,
$$\text{Tr}_{\Fp_q/\Fp_p}(\displaystyle\sum_{i \neq \pm 1}a_i(x_i+x_{-1}b_{-i}-y_i)+a_{-1}(x_{-1}-y_{-1})) = 0\text{ for all }a_2, \ldots, a_{-1} \in \Fp_q.$$ Similarly to case 3 above, this is equivalent to $y_{-1}=x_{-1}$ and $y_i = x_i+x_{-1}b_{-i}$ for all $i, i \neq \pm 1.$ Now let put $a_1 \neq 0$ and $a_i = 0$ for all $i, i \neq 1$. Then 
\begin{equation*}
\text{Tr}_{\Fp_q/\Fp_p}(y_1a_1) = \text{Tr}_{\Fp_q/\Fp_p}(x_1a_1 + \displaystyle\sum_{i \neq \pm 1}x_i(-a_1b_i)+x_{-1}(-a_1)\displaystyle\sum b_ib_{-i})
\end{equation*}
for all $a_1 \in \Fp_q$, which is equivalent to the condition $$y_1 = x_1 - \displaystyle\sum_{i \neq \pm 1}x_ib_i-x_{-1}\displaystyle\sum b_ib_{-i}.$$
It follows from the linearity of the map $\mathrm{Tr}$ that the relations on $x$ and $y$ above are not only necessary conditions for the equality $b^{-1}.\psi_x = \psi_y$, but are also sufficient ones.

In cases (i) and (ii), $x_{-1} = 0$ and $y_i = x_i$ for all $i\neq 1$, and also $y_1 = x_1 - \displaystyle\sum_{i \neq \pm 1}b_ix_i.$ That is, each orbit has form $(*, x_2, \ldots, x_{-2},0)$, and case (i) gives $q^{k-2}$ different irreducible characters of dimension~1, case (ii) gives $q^{k-2}\dfrac{q^{k-2} - 1}{q - 1}$ irreducible characters of dimension $q$. In case (iii), the scalars $b_i$ are uniquely determined by the respective $y_i,~~i \neq \pm 1$, while $y_1$ and $y_{-1}$ are uniquely determined by $x$ and $b$ due to the relations on $x$ and $y$ above. Thus, there are $(q-1)q$ different irreducible characters of dimension $[\Hei_{2}(B_{n+1}):A]=q^{2n-1}$ (and of dimension $[\Hei_{2}(D_{n+1}):A]=q^{2n-2}$ respectively).

So, by the Kirillov's bijection we get the follows.
\[
\renewcommand{\arraystretch}{3}
\begin{tabular}{|c|c|c|}
    \hline
    The dimension &  \multicolumn{2}{|c|}{\text{The number of coadjoint orbits}} \\
    \cline{2-3}
    & \text{Type $B_n$} & \text{Type $D_n$} \\
    \hline
    0 & $q^{2n-3}$ & $q^{2n-4}$ \\
    \hline
    1 & $q^{2n-3}\dfrac{q^{2n-3} - 1}{q - 1}$ & $q^{2n-4}\dfrac{q^{2n-4} - 1}{q - 1}$ \\
    \hline
    $2n-3$ & $(q-1)q$ &  0\\
    \hline
    $2n-4$ & 0 & $(q-1)q$ \\
    \hline
\end{tabular}
\]

Thus, we have 
\theobbp{The\label{theo:2_Hei_Isaacs} Isaacs\textup' conjecture holds for orthogonal $2$-layered Heisenberg groups over a finite fields of characteristic greater than $3$.}

\bigskip\textsc{Mikhail Ignatev: National Research University Higher School of Economics,\break\indent Pokrovsky Boulevard 11, 109028, Moscow, Russia}

\emph{E-mail address}: \texttt{mihail.ignatev@gmail.com}

\medskip\textsc{Leonid Titov: National Research University Higher School of Economics,\break\indent Pokrovsky Boulevard 11, 109028, Moscow, Russia}

\emph{E-mail address}: \texttt{titovleonid1986@gmail.com}


\begin{thebibliography}{XXXXXX}\addcontentsline{toc}{subsection}{References}

%\bibitem[An95]{Andre95} C.A.M. Andr\'e. Basic sums of coadjoint orbits of the unitriangular group. J.~Algebra \textbf{176} (1995), 959--1000.

%\bibitem[An01]{Andre01} C.A.M. Andr\'e. The basic character table of the unitriangular group. J.~Algebra \textbf{241} (2001), 437--471.

\bibitem[AN06]{AndreNeto06} C.A.M Andr\'e, A.M. Neto. Super-characters of finite unipotent groups of types $B_n$, $C_n$ and~$D_n$. J. Algebra \textbf{305} (2006), 394--429.



%\bibitem[Ba99]{Baranov99} A. Baranov. Finitary simple Lie algebras, J. Algebra \textbf{219} (1999), 299--329.

\bibitem[Bou02]{Bou} N. Bourbaki. Lie groups and Lie algebras. Chapters 4--6, Springer, 2002.


\bibitem[Cu24]{Cushman24} R.A. Cushman. Momentum map for the Heisenberg group. Symmetry \textbf{16} (2024), no. 8, article no. 1054.


\bibitem[De24]{Deloup24} F.L. Deloup. On Heisenberg groups, arXiv: \texttt{math.GR/2409.03399}.


\bibitem[DM08]{DarafshehMisaghian08} M.R. Darafsheh, M. Misaghian. On the ordinary irreducible characters of the Heisenberg group and a similar special group. Algebra Colloquium \textbf{15} (2008), no, 3, 471--478.




%\bibitem[BGR]{BorhoGabrielRentschler} W. Borho, P. Gabriel, R. Rentschler. Primideale in Einh$\ddot{\mathrm{u}}$llenden aufl$\ddot{\mathrm{o}}$sbarer Lie-Algebren, Lecture Notes in Math. \textbf{357}. Springer--Verlag, Berlin, 1973.


%\bibitem[Bou]{Bou} N. Bourbaki. Algebra II. Chapters 4--7. Elements of mathematics (Berlin). Springer--Verlag, Berlin, 2003.


%\bibitem[DP99]{DimitrovPenkov99} I. Dimitrov, I. Penkov. Weight modules of direct limit Lie algebras, Int. Math. Res. Notes \textbf{5} (1999), 223--249.


%\bibitem[Di96]{Dixmier1} J. Dixmier. Enveloping algebras. Grad. Stud. in Math. \textbf{11}. AMS, 1996.


\bibitem[Ev10]{Evseev10} A. Evseev. Reduction for characters of finite algebra groups. J. Algebra \textbf{325} (2010),\break 321--351.

%\bibitem[GAP]{GAP} Groups Algorithms Programming, software, \texttt{https://www.gap-system.org/install/}.

%\bibitem[Gi11]{Ginzburg11} D. Ginzburg. On the representations of $q$ unipotent groups over a finite field. Israel J. Math. \textbf{181} (2011), 387--422.

\bibitem[GMR16]{GoodwinMoschRoehrle16} S.M. Goodwin, P. Mosch, G. R\"ohrle. On the coadjoint orbits of maximal unipotent subgroups of reductive groups. Transformation Groups \textbf{21} (2016), 399--426.
%GoodwinTables
\bibitem[GLMP16]{GoodwinLeMagaardPaolini16} S.M. Goodwin, T. Le, K. Magaard, A. Paolini. Constructing characters of Sylow\break $p$-subgroups of finite Chevalley groups. J. Algebra \textbf{468} (2016), 395--439.
%GoodwinIdea
\bibitem[GLM17]{GoodwinLeMagaard} S.M. Goodwin, T. Le, K. Magaard. The generic character table of a Sylow $p$-subgroup of a finite Chevalley group of type $D_4$. Commun.
Algebra \textbf{45} (2012), 5158--5179.
%D4

\bibitem[Hi60]{Higman60} G. Higman. Enumerating $p$-groups. I. Inequalities. Proc. London Math. Soc. (3) \textbf{10} (1960), 24--30.

%\bibitem[Ig09]{Ignatev09} M.V. Ignat'ev, {\it Orthogonal subsets of classical root systems and coadjoint orbits of unipotent groups}, Math. Notes \textbf{86} (2009), no. 1, 65--80, see also arXiv: \texttt{math.RT/0904.2841}.

%\bibitem[Ig11]{Ignatev11} M.V. Ignat'ev, {\it Orthogonal subsets of root systems and the orbit method}, St. Petersburg Math.~J., \textbf{22} (2011), no. 5, 777--794, see also arXiv: \texttt{math.RT/1007.5220}.

%\bibitem[Ig12]{Ignatyev12} M.V. Ignatyev. Combinatorics of $B$-orbits and Bruhat--Chevalley order on involutions. Transformation Groups \textbf{17} (2012), no. 3, 747--780; arXiv: \texttt{math.RT/1101.2189}.

%\bibitem[Ig19]{Ignatyev19} M.V. Ignatyev. Centrally generated ideals of $U(\nt)$ in types $B$ and $D$. Transformation Groups \textbf{24} (2019), no. 4, 1067--1093; arXiv: \texttt{math.RT/1907.04219}.

%\bibitem[Ig09]{Ignatev09} M.V. Ignatev. Subregular characters of the unitriangular group over a finite field. J. Math. Sci. \textbf{156} (2009), no. 2, 276--291.


\bibitem[IPa09]{IgnatevPanov09} M.V. Ignatev, A.N. Panov. Coadjoint orbits of the group $\mathrm{UT}(7, K)$. J. Math. Sci. \textbf{156} (2009), no. 2, 292--312.

%\bibitem[IPe16]{IgnatyevPenkov16} M.V. Ignatyev, I. Penkov. Infinite Kostant cascades and centrally generated primitive ideals of $U(\nt)$  in types $A_{\infty}$, $C_{\infty}$. J. Algebra \textbf{447} (2016), 109--134; arXiv: \texttt{math.RT/1502.05486}.

\bibitem[IPe25]{IgnatevPetukhov25} M. Ignatev, A. Petukhov. Coadjoint orbits of low dimension for nilradicals of Borel subalgebras in classical types. Expositiones Mathematicae, submitted, arXiv: \texttt{math.RT/2507.20332}. 

%\bibitem[IS23]{Sur23} M.V. Ignatev, M.A. Surkov. Rook placements in $G_2$ and $F_4$ and  and associated coadjoint orbits. Commun. Math. \textbf{30} (2022), no. 2, 129--149.

%\bibitem[IS21]{IgnatyevShevchenko21} M.V. Ignatyev, A.A. Shevchenko. Centrally generated primitive ideals of $U(\nt)$ for exceptional types. J. Algebra \textbf{565} (2021), 627--650.

\bibitem[Is07]{Isaacs07} I.M. Isaacs. Counting characters of upper triangular groups. J. Algebra \textbf{315} (2007),\break 698--719.


%\bibitem[Ka77]{Kazhdan77} D. Kazhdan. Proof of Springer's hypothesis. Israel J. Math. \textbf{28} (1977), 272--286

\bibitem[Ki62]{Kirillov62} A.A. Kirillov. Unitary representations of nilpotent Lie groups. Russian Math. Surveys \textbf{17}~(1962), 53--110.

%\bibitem[Ki04]{Kirillov04} A.A. Kirillov. Lectures on the orbit method. Grad. Stud. in Math. \textbf{64}, AMS, 2004.


%\bibitem[Kr00]{Kraft00} H. Kraft. Geometric methods in invariant theory. IO NFMI, 2000.

%\bibitem[Ko12]{Kostant12} B. Kostant. The cascade of orthogonal roots and the coadjoint structure of the nilradical of a Borel subgroup of a semisimple Lie group, Moscow Math. J. \textbf{12} (2012), no. 3, 605--620.

%\bibitem[Ko13]{Kostant13} B. Kostant. Center of $U(\nt)$, cascade of orthogonal roots and a construction of Lipsman--Wolf. In: Lie groups: structure, actions and representations, Progr. in Math. \textbf{306}. Birkh$\ddot{\mathrm{a}}$user, 2013, 163--174.
 


\bibitem[Le10]{Le10} T. Le. Counting irreducible representations of large degree of the upper triangular groups. J. Algebra \textbf{324} (2010), 1803--1817.

\bibitem[LMP20]{LeMagaardPaolini20} T. Le, K. Magaard, A. Paolini. On the characters of Sylow $p$-subgroups of finite Chevalley groups $G(p^f)$ for arbitrary primes. Mathematics of Computation \textbf{89} (2020), 1501--1524.


\bibitem[Le74]{Lehrer74} G.I. Lehrer. Discrete series and the unipotent subgroup. Compositio Math., \textbf{28} (1974),\break fasc. 1, 9--19.

\bibitem[LN09]{LidlNiederreiter09} R. Lidl, H. Niederreiter. Finite fields. 
Encyclopedia of Mathematics and its Applications \textbf{20}, Cambridge University Press, Cambridge, 2009.

\bibitem[Lou11]{Loukaki11} M. Loukaki. Counting characters of small degree in upper unitriangular groups. J. Pure Appl. Algebra \textbf{215} (2011), no. 2, 154--160.

%\bibitem[Mu05]{Mukherjee05} S. Mukherjee. Coadjoint orbits for $A_{n-1}^+$, $D_n^+$ and $D_n^+$, arXiv: \texttt{math.RT/0501332}.

\bibitem[Mb11]{Marberg11} E. Marberg. Combinatorial methods of character enumeration for the unitriangular group. J. Algebra \textbf{345} (2011), 295--323.

\bibitem[Mi11]{Misaghian11} M. Misaghian. The representations of the Heisenberg group over a finite field. Armenian Journal of Mathematics \textbf{3} (2011), no. 4, 162--173.

\bibitem[Mj97]{Marjoram97} M. Marjoram. Irreducible characters of a Sylow $p$-subgroups of the orthogonal group. Extracta Mathematicae \textbf{12} (1997), no. 3, 315--319.

\bibitem[Mj97']{Marjoram97'} M. Marjoram. Irreducible characters of Sylow p-subgroups of classical groups. PhD thesis,
National University of Ireland, Dublin, 1997.

\bibitem[Mj99]{Marjoram99} M. Marjoram. Irreducible characters of small degree of the unitriangular group. Irish Math. Soc. Bull. \textbf{42} (1999), 21--31.


\bibitem[Ni20]{Nien20} C. Nien. Characters of 2-layered Heisenberg groups. Linear and Multilinear Algebra,\break \texttt{doi}: \texttt{10.1080/03081087.2020.1849005} (2020).

\bibitem[PS15]{PakSoffer15} I. Pak, A. Soffer. On Higman's $k(U_n(q))$ conjecture, arXiv: \texttt{math.CO/1507.00411}.


%\bibitem[Pa08]{Panov08} A.N. Panov. Involutions in $S_n$ and associated coadjoint orbits. J. Math. Sci. \textbf{151} (2008), 3018--3031.

%\bibitem[Pa09]{Panov09} A.N. Panov. On the index of certain nilpotent Lie algebras. J. Math. Sci. \textbf{161} (2009),\break no. 1, 122--129.

%\bibitem[Sa09]{Sangroniz09} J. Sangroniz. Irreducible characters of large degree of Sylow $p$-subgroups
of classical groups. J. Algebra \textbf{321} (2009), 1480--1496.

\bibitem[Se77]{Serre77} J.-P. Serre. Linear represenations of finite groups. Grad. Texts in Math. \textbf{42}, Springer--Verlag, New York, 1977.


\bibitem[Sur26]{Sur26} M. A. Surkov. Classification of coadjoint orbits for the maximal unipotent subgroup in the simple group of type $F_4$. Commun. Math. \textbf{34} (2026), no. 1, article no. 4.

\bibitem[Sz06]{Szechtman06} F. Szechtman. Irreducible characters of Sylow subgroups of symplectic and unitary groups. J. Algebra \textbf{303} (2006), 722--730.

\bibitem[Ve26]{Venchakov26} M. Venchakov. Rook placements and coadjoint orbits for maximal unipotent subgroups of finite symplectic groups, arXiv: \texttt{math.RT/2602.14933}.

\bibitem[VLA03]{VeraLopezArregi03} A. Vera-L\'opez, J.M. Arregi. Conjugacy classes in unitriangular matrices. Linear Algebra and its Applications \textbf{370} (2003), 85--124.

%\bibitem[Ve70]{Vergne} M. Vergne. Construction de sous-alg\`{e}bres subordonn\'{e}es \`{a} un \'{e}l\'{e}ment du~dual d'une alg\`{e}bre de Lie r\'{e}soluble. C. R. Acad. Sci. Paris Ser. A--B \textbf{270} (1970), A173--A175.


\end{thebibliography}
\end{document}